\documentclass[a4paper,12pt]{article}
\usepackage[utf8]{inputenc}

\usepackage{amssymb}
\usepackage{amsthm}
\usepackage{amsmath}
\usepackage{hyperref}
\usepackage[open,openlevel=1]{bookmark}
\newtheorem{theorem}{Theorem}[section]

\newtheorem{corollary}{Corollary}[section]

\newtheorem{remark}[theorem]{Remark}
\newtheorem{conjecture}{Conjecture}

\newtheorem*{acknowledgement}{Acknowledgement}
\theoremstyle{definition}

\begin{document}

\title{A Systolic Inequality for the Filling Area Conjecture}

\author{Chaitanya Ambi\footnote{Chennai Mathematical Institute, H1, SIPCOT IT Park, Siruseri,
Kelambakkam 603103. India. Email: chaitanya.ambi@gmail.com}}
\maketitle

\begin{abstract}
\indent

   We prove an upper bound on the systolic ratio of an orientable isometric filling of the circle equipped with a Riemannian metric. The bound depends only on the genus of isometric filling. We also apply the bound to the class of orientable isometric filling with a certain lower bound on the systole. We deduce that the Filling Area Conjecture holds true for this class when the genus is sufficiently large.
\end{abstract}

\section{Introduction.}\label{intro}
\indent

    The \textit{filling volume} of a compact, orientable $n$-dimensional manifold $N $ with a piecewise smooth Riemannian metric $\mathcal{H}$ is an important metric invariant discovered by Gromov in \cite{gromov}. It is defined for $n\geq 2$ as
    \begin{align}\label{FillVol_def}
 FillVol(N,\mathcal{H}):= \inf_{(M,\mathcal{G})}\big \lbrace vol_{n+1}(M,\mathcal{G}) : N=\partial M \text{ and } (dist_{\mathcal{G}})|_{\partial M}=dist_{\mathcal{H}} \big \rbrace,
\end{align}
    where $dist_{\mathcal{G} }$ and $dist_{\mathcal{H} }$ denote the distances induced by $\mathcal{G}$ and $\mathcal{H}$, respectively.
    Whenever \begin{equation*}
              (dist_{\mathcal{G}})|_{\partial M}=dist_{\mathcal{H}},
             \end{equation*}
    we say that $M$ fills $N$ \textit{isometrically}. We refer the reader to \cite{katz} for more details.\\

    For $n=1$, Gromov's Filling Area Conjecture (see \cite{gromov}) may be stated as follows.
    \begin{conjecture}\label{conj-gromov}
    For a smooth Riemannian metric $\mathcal{H}$ on the circle $\mathbb{S}^{1}$,  we have
    \begin{equation*}
     FillVol(\mathbb{S}^{1},\mathcal{H}) =\dfrac{vol_{1}(\mathbb{S}^{1},\mathcal{H})^{2}}{2\pi}.
    \end{equation*}
    \end{conjecture}
    If one assumes that a minimal filling is diffeomorphic to a hemisphere, then Conj. (\ref{conj-gromov}) can be deduced from Pu's systolic inequality (see \cite{katz}, Eq. (5.2.2), p. 40). Note the upper bound
   \begin{align}\label{eq_const}
    FillVol(\mathbb{S}^{1},\mathcal{H})\leq \dfrac{vol_{1}(\mathbb{S}^{1},\mathcal{H})^{2}}{2\pi},
   \end{align}
   obtained by considering the isometric filling by the hemisphere of constant curvature. In \cite{bangert}, it is proven that Conj. (\ref{conj-gromov}) holds true for all genus one fillings. (See also \cite{ivanov} for the Finslerian case.)\\

       We shall denote by $\Sigma_{g}$ a compact surface of genus $g\geq 1$ throughout this article. By virtue of the upper bound in Eq. (\ref{eq_const}), we may restrict our attention to following the class $\mathcal{I}_{g}(\mathbb{S}^{1},\mathcal{H})$ of piecewise smooth, orientable and \textit{isometric} filling $(M,\mathcal{G})$ of $(\mathbb{S}^{1},\mathcal{H})$ of fixed genus $g\geq 1$ (defined so that $\chi(M)=1-2g$).
    \begin{multline*}
     \mathcal{I}_{g}(\mathbb{S}^{1},\mathcal{H}):=\bigg\lbrace (M,\mathcal{G}) : M \text{ has genus  }g \text{ and }vol_{2}(M,\mathcal{G}) < \dfrac{vol_{1}(\mathbb{S}^{1},\mathcal{H})^{2}}{2\pi}  \bigg\rbrace.
    \end{multline*}
    Thus, to prove Conj. (\ref{conj-gromov}) for orientable isometric fillings is tantamount to showing that the class $\mathcal{I}_{g}(\mathbb{S}^{1},\mathcal{H})$ is empty for each $g\geq 1$.\\

    Consider a compact surface (with or without boundary) carrying a Riemannian metric $(\Sigma, \mathcal{G})$. Whenever $\Sigma$ is not simply connected, we define a \textit{systole} to be a closed curve with the shortest possible length which is not null-homotopic. We denote its length by $sys\pi_{1}(\Sigma, \mathcal{G}) $. The \textit{systolic ratio} of $(\Sigma, \mathcal{G})$ (denoted by $SR(\Sigma, \mathcal{G})$) is defined as follows.
    \begin{equation*}
     SR(\Sigma, \mathcal{G}):= \dfrac{(sys\pi_{1}(\Sigma, \mathcal{G}))^{2}}{vol_{2}(\Sigma, \mathcal{G})}.
    \end{equation*}
    When $\partial \Sigma = \emptyset$, we further denote the supremum of $SR(\Sigma, \mathcal{G})$ over all metrics $\mathcal{G}$ by $SR(\Sigma)$. For estimates of the systolic ratio of surfaces without boundary, we refer the reader to \cite{katz}.\\

    As mentioned before, Pu's systolic inequality yields the bound given by Eq. (\ref{eq_const}); see \cite{gromov}. This suggests intimate connections between systolic and filling inequalitites. Indeed, we prove the following.
   \begin{theorem}\label{main}
    Let $\mathcal{H}$ be a smooth Riemannian metric on $\mathbb{S}^{1}$. Consider an orientable and piecewise smooth isometric filling $(M,\mathcal{G})$ of $(\mathbb{S}^{1},\mathcal{H})$. Assume that $\pi_{1}(M)$ is non-trivial and put $g:=(1-\chi(M))/2$. Then the systolic ratio of $(M,\mathcal{G})$ satisfies
    \begin{equation}\label{eq_main}
    SR(M,\mathcal{G}) \leq SR(\Sigma_{g+1})/2.
    \end{equation}
   \end{theorem}
   \indent

   The above result applies, in particular, to every possible minimal filling of non-zero genus. Furthermore, Thm. (\ref{main}) allows us to infer the following.
   \begin{corollary}\label{cor}
    Let the notation and assumptions be as in Theorem \ref{main}. There exists an integer $g_{0}\geq 1$ (independently of $(\mathbb{S}^{1},\mathcal{H})$) such that whenever $g\geq g_{0}$ and the systole of $(M,\mathcal{G})$ satisfies
    \begin{equation}\label{eq_cor}
    sys\pi_{1}(M,\mathcal{G}) > \big (\dfrac{\log g}{2\pi \sqrt{g}}\big )\cdot vol_{1}(\mathbb{S}^{1},\mathcal{H}),
    \end{equation}
    we have $(M,\mathcal{G}) \not \in \mathcal{I}_{g}(\mathbb{S}^{1},\mathcal{H})$.
   \end{corollary}
   Thus, every orientable isometric filling of the circle with sufficiently large genus and systole exceeding the right hand side of Eq. (\ref{eq_cor}) satisfies Conj. \ref{conj-gromov}.
    \begin{remark}
    There exist isometric fillings with arbitrarily large genus whose area is marginally greater than the filling area. Consider, for instance, an isometric filling placed on the top of a small cylinder with totally geodesic boundary $\mathbb{S}^{1}$. If we attach $g\geq 1$ small handles to the hemisphere away from the boundary, the area of the resulting surface is close to the initial area, whereas the genus has increased by $g$.\\

    Thus, it is the exact value of filling volume of a manifold which is  important.
    \end{remark}
    \begin{remark}
    \indent

     Interestingly, we may also consider Def. (\ref{FillVol_def}) in relation with cobordism. Set
     \begin{equation}\label{eq_cobord}
      (N, \mathcal{H})= (N_{1}, \mathcal{H}_{1}) \amalg (N_{2}, \mathcal{H}_{2}),
     \end{equation}
      the disjoint union two $n$-dimensional compact, orientable and connected Riemannian manifolds which are \textit{cobordant}. One may consider the infimal volume of an $(n+1)$-dimensional Riemannian manifold filling both these manifolds \textit{isometrically}. Note that if $N_{2}=\emptyset$, one recovers Def. (\ref{FillVol_def}) for $(N_{1}, \mathcal{H}_{1})$.
    \end{remark}
    \begin{acknowledgement}
     The author would like to thank Chennai Mathematical Institute for support by a post-doctoral fellowship.
    \end{acknowledgement}
\section{Proof of Theorem \ref{main}.}
\begin{proof}\label{proof_main}
\indent

    For a orientable and piecewise smooth isometric filling $(M,\mathcal{G})$ of $(\mathbb{S}^{1},\mathcal{H})$ (with $g\geq 1$), we have $\partial M=\mathbb{S}^{1}$ and $(dist_{\mathcal{G}})|_{\partial M}=dist_{\mathcal{H}}$. This implies, in particular, that $\partial M$ is a smooth geodesic in $M$ with
    \begin{equation*}                                                                                                                                                                                                                                                                                        vol_{1}(\partial M,\mathcal{G}|_{\partial M})=vol_{1}(\mathbb{S}^{1}, \mathcal{H}).                                                                                                                                                                                                                                                                                       \end{equation*}
    Note that $\partial M$ is not null-homotopic because $\pi_{1}(M)$ is non-trivial. Also, there exists a systole in $M$ whose length satisfies $0<sys\pi_{1}(M,\mathcal{G})\leq vol_{1}(\mathbb{S}^{1}, \mathcal{H})$.\\

    We construct a closed surface $(\Sigma_{g+1},\mathcal{G}_{+})$ from $(M, \mathcal{G})$ as follows. Let $\partial M$ carry the orientation compatible with the one on $M$. Choose any two antipodal points $p$ and $p'$ on $\partial M$. Let $q$ be a point on $\partial M$ located at a distance $sys\pi_{1}(M,\mathcal{G})$ from $p$. The point antipodal to $q$ will be denoted by $q'$. As the geodesic arcs $pq$ and $p'q'$ both have length equal to $sys\pi_{1}(M,\mathcal{G})$, we may glue these together. Finally, we glue $pq'$ and $qp'$ together. The piecewise smooth metric determined by $\mathcal{G}$ on the resulting closed and orientable surface (which is homeomorphic to $\Sigma_{g+1} $) will be denoted by $ \mathcal{G}_{+}$.\\

    The closed curve in $(\Sigma_{g+1}, \mathcal{G}_{+}) $ determined by $pp'$ with length $sys\pi_{1}(M,\mathcal{G})$ must be a geodesic. (This holds true because $(M,\mathcal{G})$ was an isometric filling to begin with.) Moreover, since $\pi_{1}(M)$ is non-trivial, this closed geodesic must be a \textit{systole} in $(\Sigma_{g+1}, \mathcal{G}_{+})$. It is also clear that
    \begin{equation*}
    vol_{2}(\Sigma_{g+1}, \mathcal{G}_{+})=vol_{2}(M,\mathcal{G}).
    \end{equation*}
    Hence, we have
    \begin{equation}\label{eq_SR+}
     SR(\Sigma_{g+1}, \mathcal{G}_{+})= SR(M, \mathcal{G}).
    \end{equation}
      Next, we construct a closed and \textit{non-orientable} surface $(\Sigma,\mathcal{G}_{-})$ from $(M, \mathcal{G})$ by gluing $pq$ with $q'p'$ and $pq'$ with $p'q$ (in that order). Observe that the orientable double cover of $(\Sigma,\mathcal{G}_{-})$ has genus $(g+1)$. Moreover, it carries the unique piecewise smooth metric $\tilde{\mathcal{G}}_{-}$ which makes the map
      \begin{equation*}
       (\Sigma_{g+1},\tilde{\mathcal{G}}_{-})\twoheadrightarrow (M, \mathcal{G})
      \end{equation*}
      a Riemannian covering. Arguing exactly as in the case of $(\Sigma_{g+1}, \mathcal{G}_{+})$, we get
       \begin{equation}\label{eq_SR-}
     SR(\Sigma, \mathcal{G}_{-})= SR(M, \mathcal{G}).
    \end{equation}
    But now,
    \begin{equation*}
       sys\pi_{1}(\Sigma_{g+1},\tilde{\mathcal{G}}_{-})\geq 2 sys\pi_{1}(\Sigma,\mathcal{G}_{-}).
    \end{equation*}
    By Eq. (\ref{eq_SR-}), we get
    \begin{align*}
     SR(\Sigma_{g+1},\tilde{\mathcal{G}}_{-})\geq \dfrac{4(sys\pi_{1}(\Sigma,\mathcal{G}_{-}))^{2}}{2 vol_{2}(\Sigma,\mathcal{G}_{-})},\\
     =2 SR(M, \mathcal{G}).
    \end{align*}
    In view of Eq. (\ref{eq_SR+}), we see that
     \begin{equation*}
       SR(\Sigma_{g+1},\tilde{\mathcal{G}}_{-})\geq 2 SR(\Sigma_{g+1}, \mathcal{G}_{+}).
    \end{equation*}
    Considering the definition of $SR(\Sigma_{g+1})$ together with Eq. (\ref{eq_SR+}), we obtain
    \begin{equation*}
     SR(\Sigma_{g+1})\geq 2 SR(M, \mathcal{G}).
    \end{equation*}
    This proves Eq. (\ref{eq_main}) as required.
\end{proof}
\section{Proof of Corollary \ref{cor}.}
\indent

   We maintain the notation as in Section \ref{intro}.
\begin{proof}
\indent

   Assume the contrary that there are isometric fillings in the class $\mathcal{I}_{g}(\mathbb{S}^{1},\mathcal{H})$ for arbitrarily large $g\geq 1$ which satisfy Eq. (\ref{eq_cor}). For such an integer $g$, let $(M, \mathcal{G})\in \mathcal{I}_{g}(\mathbb{S}^{1},\mathcal{H})$ be any piecewise smooth, orientable and isometric filling of the circle $(\mathbb{S}^{1},\mathcal{H})$. Thus,
   \begin{equation*}
    vol_{2}(M, \mathcal{G})< (vol_{1}(\mathbb{S}^{1},\mathcal{H}))^{2}/2\pi.
   \end{equation*}
   We claim that $\partial M$ is \textit{not} a systole in $(M, \mathcal{G})$. Indeed, if it were so, then $sys\pi_{1}(M, \mathcal{G})=vol_{1}(\mathbb{S}^{1},\mathcal{H})$. Hence, the preceeding equation would force that
   \begin{align*}
    SR(M, \mathcal{G})> \dfrac{(vol_{1}(\mathbb{S}^{1},\mathcal{H}))^{2}}{(vol_{1}(\mathbb{S}^{1},\mathcal{H}))^{2}/2\pi},\\
    = 2\pi.
   \end{align*}
   Combined with Thm. \ref{main}, we would get
   \begin{equation*}
    2\pi < SR(M, \mathcal{G})\leq SR(\Sigma_{g+1})/2.
   \end{equation*}
   This stark contradiction with the basic estimate $SR(\Sigma_{g+1})\leq 4/3$ proves the claim (see \cite{katz}, Eq. (11.2.1)).\\

   Next, we glue along $\mathbb{S}^{1}\approx \partial M$ the (totally geodesic) boundary of hemisphere of constant sectional curvature equal to $(2\pi/vol_{1}(\mathbb{S}^{1},\mathcal{H}))^{2}$. We denote the piecewise smooth metric on the resulting closed, orientable surface (which is homeomorphic to $\Sigma_{g}$) by $\mathcal{G}_{0}$. Note that
   \begin{equation*}
    vol_{2}(\Sigma_{g}, \mathcal{G}_{0})= vol_{2}(M, \mathcal{G})+ (vol_{1}(\mathbb{S}^{1},\mathcal{H}))^{2}/2\pi.
   \end{equation*}
   Hence,
    \begin{equation}\label{eq_vol}
    2\cdot vol_{2}(M, \mathcal{G})\leq vol_{2}(\Sigma_{g}, \mathcal{G}_{0})\leq (vol_{1}(\mathbb{S}^{1},\mathcal{H}))^{2}/\pi.
   \end{equation}
   As $\partial M$ is not a systole in $M$, we have
   \begin{equation*}
    sys\pi_{1}(\Sigma_{g}, \mathcal{G}_{0})= sys\pi_{1}(M, \mathcal{G}).
   \end{equation*}
   If we square the above equation and use Eq. (\ref{eq_vol}), we obtain
   \begin{multline*}
    \dfrac{\pi \cdot (sys\pi_{1}(M, \mathcal{G}))^{2}}{vol_{1}(\mathbb{S}^{1},\mathcal{H})^{2}} \leq SR(\Sigma_{g}, \mathcal{G}_{0})\leq \dfrac{SR(M, \mathcal{G})}{2}\leq \dfrac{SR(\Sigma_{g+1})}{4}.
   \end{multline*}
   Here, we have used Thm. \ref{main} for the rightmost inequality. In particular,
   \begin{equation*}
   \bigg (\dfrac{sys\pi_{1}(M, \mathcal{G})}{vol_{1}(\mathbb{S}^{1},\mathcal{H})}\bigg )^{2}\leq \dfrac{SR(\Sigma_{g+1})}{4\pi}.
   \end{equation*}
   In view of the lower bound given by Eq. (\ref{eq_cor}), this implies
   \begin{equation}\label{eq_final}
    \frac{1}{\pi} < \dfrac{g\cdot SR(\Sigma_{g+1})}{\log^{2}(g)}.
   \end{equation}
   But now, we also have (see \cite{katz}, Thm. 11.3.1, Eq. (11.3.1), p. 88)
    \begin{equation}\label{eq_sys}
    \limsup_{g\rightarrow\infty}\dfrac{(g+1)\cdot SR(\Sigma_{g+1})}{\log^{2} (g+1)}\leq \frac{1}{\pi}.
    \end{equation}
    In view of the above Eq. (\ref{eq_sys}), we see that Eq. (\ref{eq_final}) cannot hold for arbitrarily large value of $g$. This contradiction shows that there must be an integer $g_{0}\in \mathbb{N}$ such that for $g\geq g_{0}$, Cor. (\ref{cor}) holds true.
\end{proof}
\bibliographystyle{plain}
\bibliography{Filling-Area}
\end{document}